\documentclass[12pt,reqno]{amsart}
\usepackage{amssymb,verbatim,enumerate,ifthen}
\usepackage[mathscr]{eucal}
\usepackage[utf8]{inputenc}
\usepackage[T1]{fontenc}
\def\R{\mathbb{R}}
\def\Q{\mathbb{Q}}

\def\sign{\mathop{\mbox{\rm sign}}\nolimits}

\newtheorem{theorem}{Theorem}
\newtheorem*{theorem*}{Theorem}
\def\Thm#1#2{\ifthenelse{\equal{#1}{*}}{\begin{theorem*}#2\end{theorem*}}
             {\begin{theorem}\label{T#1}#2\end{theorem}}}
\newtheorem{Atheorem}{Theorem}

\def\thm#1{Theorem~\ref{T#1}}

\newtheorem{proposition}[theorem]{Proposition}
\newtheorem*{proposition*}{Proposition}
\def\Prp#1#2{\ifthenelse{\equal{#1}{*}}{\begin{proposition*}#2\end{proposition*}}
             {\begin{proposition}\label{P#1}#2\end{proposition}}}

\newtheorem{corollary}[theorem]{Corollary}
\newtheorem*{corollary*}{Corollary}
\def\Cor#1#2{\ifthenelse{\equal{#1}{*}}{\begin{corollary*}#2\end{corollary*}}
             {\begin{corollary}\label{C#1}#2\end{corollary}}}

\newtheorem{lemma}[theorem]{Lemma}
\newtheorem*{lemma*}{Lemma}
\def\Lem#1#2{\ifthenelse{\equal{#1}{*}}{\begin{lemma*}#2\end{lemma*}}
             {\begin{lemma}\label{L#1}#2\end{lemma}}}

\newtheorem{remark}[theorem]{Remark}
\newtheorem*{remark*}{Remark}
\def\Rem#1#2{\ifthenelse{\equal{#1}{*}}{\begin{remark*}\rm #2\end{remark*}}
             {\begin{remark}\label{R#1}\rm #2\end{remark}}}

\newtheorem{example}[theorem]{Example}
\newtheorem*{example*}{Example}
\def\Exa#1#2{\ifthenelse{\equal{#1}{*}}{\begin{example*}\rm #2\end{example*}}
             {\begin{example}\label{Ex#1}\rm #2\end{example}}}

\def\eq#1{{\rm(\ref{E#1})}}
\def\Eq#1#2{\ifthenelse{\equal{#1}{*}}
  {\begin{equation*}\begin{aligned}[]#2\end{aligned}\end{equation*}}
  {\begin{equation}\begin{aligned}[]\label{E#1}#2\end{aligned}\end{equation}}}

\begin{document}
\begin{flushright}
\textit{Submitted to: Aequationes Math.}
\end{flushright}
\vspace{5mm}

\date{\today}

\title[Convexity with respect to families of means]{Convexity with respect to families of means}

\author[Gy. Maksa]{Gyula Maksa}
\author[Zs. P\'ales]{Zsolt P\'ales}
\address{Institute of Mathematics, University of Debrecen,
H-4010 Debrecen, Pf.\ 12, Hungary}
\email{\{maksa,pales\}@science.unideb.hu}

\subjclass[2010]{Primary 39B62; Secondary 26D07, 26D15}
\keywords{}

\dedicatory{Dedicated to the 90th birthday of Professor János Aczél}

\thanks{The research of the second author was realized in the frames of T\'AMOP 4.2.4. A/2-11-1-2012-0001
”National Excellence Program -- Elaborating and operating an inland student and researcher personal
support system”. This project was subsidized by the European Union and co-financed by the European
Social Fund. The research of both authors was also supported by the Hungarian Scientific
Research Fund (OTKA) Grant NK 81402.}

\begin{abstract}
In this paper we investigate continuity properties of functions $f:\R_+\to\R_+$ that satisfy
the $(p,q)$-Jensen convexity inequality
\Eq{*}{
  f\big(H_p(x,y)\big)\leq H_q(f(x),f(y)) \qquad(x,y>0),
}
where $H_p$ stands for the $p$th power (or Hölder) mean. One of the main results shows that
there exist discontinuous multiplicative functions that are $(p,p)$-Jensen convex for all positive
rational number $p$. A counterpart of this result states that if $f$ is $(p,p)$-Jensen convex for 
all $p\in P\subseteq\R_+$, where $P$ is a set of positive Lebesgue measure, then $f$ must be continuous.
\end{abstract}

\maketitle

\section{Introduction}

Given two intervals $I$ and $J$ in $\R$ and two two-variable means $M:I^2\to I$ and
$N:J^2\to J$, a real function $f:I\to J$ is called \emph{$(M,N)$-convex} if
\Eq{*}{
  f(M(x,y))\leq N(f(x),f(y)) \qquad(x,y\in I).
}
In the particular case when $M$ and $N$ are the arithmetic means, $(M,N)$-convexity
is termed \emph{Jensen convexity}. Jensen convex functions, in general, are not necessarily
continuous. In view of the existence of a Hamel base of $\R$ considered as a vector space
over $\Q$, one can construct a discontinuous additive function (\cite{Ham05}, \cite{Kuc85}),
such a function is automatically a discontinuous Jensen convex function. On the other hand,
certain weak regularity properties of Jensen convex functions ensure that they are also convex
and henceforth continuous. Jensen \cite{Jen05}, \cite{Jen06} verified that continuous Jensen
convex functions are convex. Bernstein and Doetsch \cite{BerDoe15} showed that upper boundedness
on a nonempty open set and Jensen convexity imply convexity. Sierpiński \cite{Sie20} proved that
upper boundedness on a set of positive Lebesgue measure and Jensen convexity are also sufficient
conditions for convexity. In the context of $(M,N)$-convexity, Zgraja \cite{Zgr03} proved that
with sufficiently regular means $M$ and $N$, weak regularity properties of $(M,N)$ convex functions
imply their continuity.

In this paper we will consider the particular case of $(M,N)$-convexity when $I=J=\R_+:=]0,\infty[$
and the means $M$ and $N$ belong to the class of power (or Hölder) means.
We recall that, for $p\in\R$ the $p$th power mean $H_p$ on $\R_+$ is defined by
\Eq{*}{
  H_p(x,y):=
  \begin{cases}
  \left(\dfrac{x^p+y^p}{2}\right)^{\frac{1}{p}} & \mbox{if $p\neq0$}\\[2mm]
  \sqrt{xy} & \mbox{if $p=0$}
  \end{cases}
  \qquad(x,y\in\R_+)
}
(cf.\ \cite{HarLitPol34}, \cite{Hol1889}).
In the sequel $(H_p,H_q)$-convex functions will simply be called \emph{$(p,q)$-Jensen convex}.
Obviously, $(1,1)$-Jensen convexity is equivalent to Jensen convexity.

The results of this paper are motivated by the following question of Matkowski \cite{Mat92}
formulated at the 30th ISFE in 1992: Is it true that $(0,0)$-Jensen convexity and $(1,1)$-Jensen convexity
imply continuity? During the symposium, a negative answer was given by Maksa \cite{Mak92}
who constructed a discontinuous $(0,0)$-Jensen convex and $(1,1)$-Jensen convex function.

The aim of this paper is to investigate continuity properties of functions that are
$(p,q)$-Jensen convex for $(p,q)\in\Pi$, where $\Pi\subseteq\R^2$ is a nonempty set.
Such functions will also be termed $\Pi$-Jensen convex.

\section{Preliminary observations}

We start with a characterization of $(p,q)$-Jensen convexity. To formulate this result,
for $(p,q)\in\R^2$ and for $f:I\subseteq\R_+\to\R_+$, define the function $f_{p,q}:I_p\to\R$ by
\Eq{*}{
  f_{p,q}(x):=
  \begin{cases}
   \sign(q)\Big(f\big(x^{\frac1p}\big)\Big)^q\qquad&\mbox{if $p\neq0$, $q\neq0$} \\[2mm]
   \sign(q)\big(f\circ\exp(x)\big)^q\qquad&\mbox{if $p=0$, $q\neq0$} \\[2mm]
   \log\circ f\big(x^{\frac1p}\big)\qquad&\mbox{if $p\neq0$, $q=0$} \\[2mm]
   \log\circ f\circ\exp(x)\qquad&\mbox{if $p=0$, $q=0$}
  \end{cases}
  \qquad (x\in I_p),
}
where
\Eq{*}{
  I_p:=
  \begin{cases}
   I^p:=\{t^p\mid t\in I\}\qquad&\mbox{if $p\neq0$}, \\[2mm]
   \log(I):=\{\log t\mid t\in I\}\qquad&\mbox{if $p=0$}.
  \end{cases}
}

The following result is very well known in the mathematical folklore. It could easily be deduced
from the results of Aczél \cite{Acz47a}.

\Thm{pq}{Let $I\subseteq \R_+$ be an interval. Then, for any $(p,q)\in\R^2$, a function
$f:I\to\R_+$ is $(p,q)$-Jensen convex if and only if $f_{p,q}$ is Jensen convex on $I_p$.}

\begin{proof}
We prove the statement only in the case $pq\neq0$. In this case, the $(p,q)$-Jensen convexity of $f$
means that
\Eq{*}{
  f\bigg(\left(\dfrac{x^p+y^p}{2}\right)^{\frac{1}{p}}\bigg)
  \leq \left(\dfrac{(f(x))^q+(f(y))^q}{2}\right)^{\frac{1}{q}} \qquad(x,y\in I).
}
Taking the $q$th power of both sides and substituting $x=u^{\frac{1}{p}}$ and $y=v^{\frac{1}{p}}$
(where $u,v\in I_p$), this inequality can be rewritten as
\Eq{*}{
  \sign(q)\left(f\bigg(\left(\dfrac{u+v}{2}\right)^{\frac{1}{p}}\bigg)\right)^q
  \leq \sign(q)\dfrac{\big(f\big(u^{\frac{1}{p}}\big)\big)^q+\big(f\big(u^{\frac{1}{p}}\big)\big)^q}{2}
   \qquad(u,v\in I_p).
}
This inequality, however, is equivalent to the Jensen convexity of the function $f_{p,q}$.

The proof in the case $pq=0$ is completely analogous. If $p=0$ then the substitutions $x=\exp(u)$
and $y=\exp(v)$ should be applied. If $q=0$, then the logarithm of the two sides of the
$(p,q)$-Jensen convexity inequality should be taken.
\end{proof}

\Cor{pq}{Let $I\subseteq \R_+$ be an open interval. If $f:I\to\R_+$ is $(p,q)$-Jensen convex
for some $q<0$ and $p\in\R$, then $f_{p,q}$ is convex on $I_p$ and, consequently $f$ is continuous,
moreover locally Lipschitz on $I$.}

\begin{proof}
If $f$ is $(p,q)$-Jensen convex, then, by the previous theorem, $f_{p,q}$ is Jensen convex on $I_p$.
On the other hand, for $q<0$, we have that $f_{p,q}\leq0$ on $I_p$. Therefore, by the
Bernstein--Doetsch theorem \cite{BerDoe15}, $f_{p,q}$ is convex. The interval $I_p$ is open, hence
$f_{p,q}$ is locally Lipschitz on $I_p$. This implies that $f$ is locally Lipschitz on $I$.
\end{proof}

In view of the above corollary, we restrict ourselves to investigate the continuity properties of
$\Pi$-Jensen convex functions, where $\Pi\subseteq \R\times[0,\infty[$ is a nonempty set.

\section{Main Results}

We recall that the real functions $a:\R\to\R$, $e:\R\to\R_+$, $\ell:\R_+\to\R$, and $m:\R_+\to\R_+$ are
called \textit{additive}, \textit{exponential}, \textit{logarithmic}, and \textit{multiplicative}
if the following Cauchy type functional equations
\Eq{*}{
  a(x+y)&=a(x)+a(y),\qquad & e(x+y)&=e(x)e(y) \qquad &(x,y&\in\R), \\
  \ell(xy)&=\ell(x)+\ell(y),\qquad & m(xy)&=m(x)m(y) \qquad &(x,y&\in\R_+)
}
are satisfied, respectively. The basic properties and of these classes of functions are
described in details in the monographs \cite{Acz66} and \cite{Kuc85}.

A functions $a:\R\to\R$ is said to be a \textit{derivation} if it is additive and satisfies
the \textit{Leibniz Rule}, that is, if $d$ fulfills the following functional equation
\Eq{*}{
  d(xy)=xd(y)+yd(x) \qquad(x,y\in\R).
}
By the Leibniz Rule $d(1)=0$ easily follows. Hence derivations vanish on $\Q$ (the rationals). Thus, a continuous derivation
is identically zero. It is a nontrivial fact that there exist nonzero derivations. Moreover, by taking an
algebraic base of $\R$, any real-valued function can uniquely be extended to a derivation (cf.\ \cite{Kuc85}).

The next theorem is one of the main results of this paper.

\Thm{M1}{For any derivation $d:\R\to\R$ and $\alpha>0$, the function
$F_{d,\alpha}:\R_+\to\R_+$ defined by
\Eq{*}{
  F_{d,\alpha}(x):=x^\alpha\exp\Big(\frac{d(x)}{x}\Big) \qquad(x\in\R_+)
}
is multiplicative and $(p,\alpha^{-1}p)$-Jensen convex for all $p\in\Q_+$(the positive rationals). Hence, if $d$ is a
\emph{nonzero derivation}, then $F_{d,\alpha}$ is a discontinuous function on $\R_+$ which is
$(p,\alpha^{-1}p)$-Jensen convex for all $p\in\Q_+$.}

\begin{proof} Observe that if $d$ satisfies the Leibniz Rule then the function $\ell(x):=d(x)/x$
is logarithmic and therefore $x\mapsto \exp(d(x)/x)$ is multiplicative. The product of two
multiplicative functions is again multiplicative, hence $F_{d,\alpha}$ is multiplicative.

The multiplicativity of $F_{d,\alpha}$ implies that $a:=\log\circ F_{d,\alpha}\circ\exp$ is additive.
Thus, by well-known properties of additive functions, $a$ is $\Q$ homogeneous, i.e., $a(px)=pa(x)$
holds for all $x\in\R$ and for all $p\in\Q$. This directly yields that
\Eq{pp}{
  F_{d,\alpha}\big(x^{p}\big)=\big(F_{d,\alpha}(x)\big)^p\qquad(x\in\R_+,\,p\in\Q).
}

First we show $F=F_{d,\alpha}$ is $(1,\alpha^{-1})$-Jensen convex. Indeed, using the inequality
between the weighted geometric mean and the weighted Hölder mean of order $\alpha^{-1}$, we get
\Eq{*}{
  F\big(H_1(x,y)\big)
  &=F\Big(\frac{x+y}{2}\Big)
   =\Big(\frac{x+y}{2}\Big)^\alpha\exp\bigg(\frac{d(x)+d(y)}{x+y}\bigg)\\
  &=\Big(\frac{x+y}{2}\Big)^\alpha
     \exp\bigg(\frac{x}{x+y}\frac{d(x)}{x}+\frac{y}{x+y}\frac{d(y)}{y}\bigg)\\
  &=\Big(\frac{x+y}{2}\Big)^\alpha
     \bigg(\frac{F(x)}{x^\alpha}\bigg)^{\frac{x}{x+y}}
         \cdot\bigg(\frac{F(y)}{y^\alpha}\bigg)^{\frac{y}{x+y}}\\
  &\leq\Big(\frac{x+y}{2}\Big)^\alpha
     \Bigg(\frac{x}{x+y}\bigg(\frac{F(x)}{x^\alpha}\bigg)^{\frac{1}{\alpha}}
         +\frac{y}{x+y}\bigg(\frac{F(y)}{y^\alpha}\bigg)^{\frac{1}{\alpha}}\Bigg)^\alpha\\
  &= \Bigg(\frac{\big(F(x)\big)^{\frac{1}{\alpha}}+\big(F(y)\big)^{\frac{1}{\alpha}}}2\Bigg)^\alpha
   = H_{\alpha^{-1}}\big((F(x),F(y)\big)\qquad(x,y\in\R_+).
}
Now, applying \eq{pp} and the $(1,\alpha^{-1})$-Jensen convexity, for $p\in\Q_+$, we obtain
\Eq{*}{
  F&(H_p(x,y))
  =F\bigg(\Big(\frac{x^p+y^p}{2}\Big)^{\frac{1}{p}}\bigg)
  =\bigg(F\Big(\frac{x^p+y^p}{2}\Big)\bigg)^{\frac{1}{p}}
  =\Big(F\big(H_1(x^p,y^p)\big)\Big)^{\frac{1}{p}} \\
  &\leq\Big(H_{\alpha^{-1}}\big(F(x^p),F(y^p)\big)\Big)^{\frac{1}{p}}
  =\Big(H_{\alpha^{-1}}\big((F(x))^p,(F(y))^p\big)\Big)^{\frac{1}{p}}
  =H_{\alpha^{-1}p}\big(F(x),F(y)\big).
}
This shows that $F_{d,\alpha}$ is $(p,\alpha^{-1}p)$-Jensen convex for all $p\in\Q_+$.
\end{proof}

\Cor{M1}{For any derivation $d:\R\to\R$ and $\alpha\geq1$,
the function $F_{d,\alpha}:\R_+\to\R_+$ defined by
\Eq{*}{
  F_{d,\alpha}(x):=x^\alpha\exp\Big(\frac{d(x)}{x}\Big) \qquad(x\in\R_+)
}
is a multiplicative and Jensen convex. Hence, if $d$ is a \emph{nonzero derivation}, then
$F_{d,\alpha}$ is a discontinuous function which is Jensen convex and multiplicative.}

In the next result we characterize those multiplicative functions that are also Jensen convex.

\Thm{M2}{A multiplicative function $m:\R_+\to\R_+$ is Jensen convex if and only if
there exists an additive function $a:\R\to\R$ such that
\Eq{mm}{
  m(t)\geq 1+a(t-1) \qquad(t\in\R_+).
}}

\begin{proof} If $m$ is a multiplicative function, then we have that $m(1)=1$.

If $m$ is Jensen convex then, as a consequence of Rod\'e's Theorem \cite{Rod78}, for every
$p>0$ there exists an additive function $a_p:\R\to\R$ such that
\Eq{*}{
  m(t)\geq m(p)+a_p(t-p) \qquad(t>0).
}
Therefore, \eq{mm} holds with $a=a_1$.

No assume that there exists an additive function $a:\R\to\R$ satisfying \eq{mm}.
Let $x,y>0$ and apply \eq{mm} for $t:=\frac{2x}{x+y}$ and $t:=\frac{2y}{x+y}$.
Adding up the inequalities so obtained side by side, we get
\Eq{*}{
  m\Big(\frac{2x}{x+y}\Big)+m\Big(\frac{2y}{x+y}\Big)
   &\geq 2 +a\Big(\frac{2x}{x+y}-1\Big)+a\Big(\frac{2y}{x+y}-1\Big)\\
   &= 2 +a\Big(\frac{2x}{x+y}+\frac{2y}{x+y}-2\Big)=2.
}
Now multiplying both sides of this inequality by $m\big(\frac{x+y}2\big)$ and using the
multiplicativity of $m$, it follows that
\Eq{*}{
  m(x)+m(y)\geq 2m\Big(\frac{x+y}2\Big),
}
which proves that $m$ is Jensen convex.
\end{proof}

The following problem is open: Find or characterize those multiplicative functions $m$
and additive functions $a$ such that \eq{mm} be valid.

\Thm{P}{Let $P\subseteq\R_+$ be a set whose internal Lebesgue measure is positive.
Let $q:P\to\R$ be an arbitrary function. If $f:\R_+\to\R_+$ is $(p,q(p))$-Jensen
convex for all $p\in P$, then $f$ is continuous (and hence, $f_{p,q(p)}$ is convex
for all $p\in P$).}

\begin{proof}
We may assume that $P$ is measurable and we may also assume that $q$ is nonnegative
because if, for some $p$, the value of $q(p)$ were negative, then \thm{pq} could be applied.
Let $0<x<y$ be fixed. Then the image $S$ of the set $P$ by the strictly increasing and
continuously differentiable mapping
\Eq{*}{
  p\mapsto H_p(x,y)
}
is of positive measure. On the other hand, in view of the $(p,q(p))$-Jensen convexity of $f$,
\Eq{*}{
  f\big(H_p(x,y)\big)\leq H_{q(p)}\big(f(x),f(y)\big)\leq\max\big(f(x),f(y)\big).
}
Therefore, $f$ is bounded from above on $S$ by $\max\big(f(x),f(y)\big)$.
Thus, for every $p\in P\setminus\{0\}$, $f_{p,q(p)}$ is bounded from above on $S^p:=\{s^{p}\mid s\in S\}$
by $\max\big(f^{q(p)}(x),f^{q(p)}(y)\big)$ if $q(p)>0$ and by $\max\big(\log(f(x)),\log(f(y))\big)$
if $q(p)=0$.

Consequently, by Sierpiński's generalization \cite{Sie20} of the Bernstein--Doetsch Theorem
\cite{BerDoe15}, for every $p\in P$, $f_{p,q(p)}$ is continuous and hence it is also convex.
\end{proof}


\end{document}